\documentclass[11pt]{article}
\usepackage{amsgen, amsmath, amsfonts, amsthm, amssymb, enumerate,amscd}

\newtheorem{defn}{Definition}[section]
\newtheorem{prop}[defn]{Proposition}

\newtheorem{thm}[defn]{Theorem}

\newtheorem{rem}[defn]{Remark}

\newcommand {\ZZ}{{\mathbb Z}}

\newcommand {\XX}{{\mathcal X}}

\newcommand {\C}{{\mathbb C}}

\newcommand {\F}{{\mathbb F}}

\newcommand {\CC}{{\mathcal C}}
\newcommand {\Q}{{\mathbb Q}}
\newcommand {\R}{{\mathbb R}}
\newcommand {\OO}{{\mathcal O}}

\newcommand {\D}{{\mathcal D}}
\newcommand {\QL}{{{\mathbb Q}_{\ell}}}
\newcommand {\bfinfty}{{\mathbf{\infty}}}

\title{Non-Archimedean   regulator maps and special values of $L$-functions.}
\author{Ramesh  Sreekantan\\University of Toronto}

\begin{document}

\maketitle
\begin{abstract}
We  define an analogue  of the  `Real' Deligne  cohomology group  at a
prime of  semi-stable  or good reduction of a variety.  We also define
regulator  maps to  this group  and formulate  a conjecture  about the
image.  This  allows us  to  formulate  a  non-Archimedian version  of
Beilinson's Hodge-$\D$-conjecture and  $S$-integral and function field
versions  of  Beilinson's global  conjectures  as  well  as a  precise
special value conjecture in the function field case. Finally we give a
few examples where these conjectures are known to be true.

\end{abstract}
\section{Introduction}

Let $X$ be a smooth projective variety over $\Q$. Beilinson
\cite{be} formulated conjectures relating special values of
$L$-functions to the $K$-theory of such varieties in terms of the
`Real' Deligne cohomology of the varieties. Roughly speaking, he
constructed a regulator map from a higher Chow group of the
variety to a real vector space, the Real Deligne cohomology, such
that the image gives a $\Q$-structure. The dimension of this
vector space is the order of pole of the Archimedean\footnote{
There seems to be some discrepancy in the literature as to the
spelling of `Archimedean' - an alternative is `Archimedian '.
However, a search on Google revealed that the former is more
popular, so that is what we have used.}
 factor of a cohomological $L$-function of the
variety. Finally the real vector space has another lattice
structure induced by the Betti and de Rham cohomology groups and
the determinant of the the change of basis matrix is related to a
special value of this $L$-function.

In this paper we define an analogue of the Deligne cohomology
group for a finite prime of good or strict semi-stable reduction.
We show there is a regulator map from the higher Chow group to
this Deligne cohomology group and show that it has similar
properties for known or conjectural reasons.

The original aim of this paper was  to formulate a version of
Beilinson's conjectures in the case of varieties over function
fields of characteristic $p$ so as  to put the results of
\cite{ko}, \cite{pal} and \cite{cs} in a general framework. Since
in this case all the primes are finite we do have such a
formulation. Further, we can formulate an $S$-integral version of
the Beilinson conjectures. Finally, we can also formulate a
precise special value conjecture in the spirit of Bloch-Kato
\cite{bk}.

While the definition of the regulator map is considerably simpler
than the Archimedean    case, in some cases there are very similar
formulas.  The link between the two seems to be through the theory
of $p$-adic uniformization.

It is generally believed that a variety should be considered to
have totally degenerate reduction at an Archimedean  place. In
particular it has semi-stable reduction, so the usual conjectures
are just the statements in the case of an Archimedean prime.

In the final section we show some examples where these conjectures
are known to be true. To a certain extent we re-interpret known
results in our terms, so we may be guilty of putting old wine in
new bottles. Further, barring the function field case, the
conjectures formulated here were perhaps implicitly, if not
explicitly, known to the experts - though as far as we are aware
they have not appeared in print - at least from this point of
view.

Finally, the correct context for the conjectures should be
motives, but we have chosen to describe them in terms of varieties
for `simplicity'.

Acknowledgements: We would like to thank Spencer Bloch for
suggesting this line of thought, Patrick Brosnan, Caterina
Consani, Najmuddin Fakhruddin, Tom Haines and Niranjan
Ramachandran for their comments and the University of Toronto for
its hospitality.

\section{The Archimedean    Case}

The usual Real Deligne cohomology $H^q_{\D}(X_{/\R},\R(q-a))$,
with $q > 2a+1$  has the following properties.

\begin{itemize}

\item 1. It is a finite dimensional real vector space with
$$dim_{\R} H^q_{\D}(X_{/\R},\R(q-a)) = -ord_{s=a} L_{\infty}(H^{q-1}(X),s)$$
the L-factor at the Archimedean    place.

\item 2. There is a regulator map $r_{\D} :
CH^{q-a}(X,q-2a)\otimes \Q \longrightarrow
H^q_{\D}(X_{/\R},\R(q-a))$

\item 3. There is a $\Q$ structure on this real vector space
induced by the Betti cohomology group $H^{q-1}_{B}(X(\C),\Q)$ and
piece of the de Rham cohomology group
$F^{q-a}H^{q-1}_{dR}(X_{/\R})$.

\item 4. The image of the regulator map is conjecturally another
$\Q$-lattice in the real vector space.

\item 5. Assuming $4$ one can compute the determinant of the
change of basis matrix with respect to these two lattices. Let
$c_{\bfinfty}(X,q,a)$ be that number. Then conjecturally
$$L^*(X,a) \sim_{\Q^*} c_{\bfinfty}(X,q,a)$$
where $L^*(X,a)$ denotes the first non-zero value of the Laurent
expansion of $L(X,s)$ at $s=a$.

\end{itemize}

We will define a $\Q$ vector space for a prime $p$ of semi-stable
or good reduction which has property $1$ owing to the work of
Consani, \cite{co}. We will then define a regulator map to this
vector space and will speculate on analogues of properties $3$,$4$
and $5$.

When $q-2a=1$ ( corresponding to $K_1$) the conjecture has to be
slightly modified - one has to add a term corresponding to the
group $B^a(X)=CH^{a}(X)/CH^{a}_{hom}(X)$ and a similar statement
holds. Our formulation takes this in to account as well.

Beilinson  \cite{be} formulated  his  conjectures in  terms of  graded
pieces  of the  $K$-theory  of the  varieties  but we  have chosen  to
formulate  it in  terms of  the  higher Chow  groups.  If  one is  not
interested in  a precise special value  conjecture it does  not make a
difference but for the exact value it could. However, in the cases for
which  we have examples  it does  not. There  is no  particularly good
reason  for our choice.   All of  the formalism  goes through  for any
candidate for motivic  cohomology, so if the higher  Chow groups fail,
it is plausible that some other candidate could succeed.

\section{Preliminaries}

Let $X$ be a smooth proper variety over a field $K$ and $\Lambda$
a discrete valuation ring with closed point $v$ and generic point
$\eta$.

By  a model $\XX$ of $X$ we mean a flat proper scheme $\XX
\rightarrow Spec(\Lambda)$ together with an isomorphism of the
generic fibre $X_{\eta}$ with $X$. Let $Y$ be the special fibre $X
\times Spec(k(v))$. We will always also make the assumption that
the model is strictly semi-stable, which means that it is a
regular model and the fibre $Y$ is a divisor with normal
crossings, the components have multiplicity one and they intersect
transversally.

We have the following picture
$$
\begin{CD}
Y @>i>> \XX @<j<< \bar{X} \\
@VVV @VVV @VVV \\
Spec(k(v)) @>>>Spec(\Lambda)@<<<Spec(\bar{\eta})
\end{CD}
$$

\section{Consani's Double Complex}

In \cite{co}, Consani defined a double complex of Chow groups of
the components of the special fibre with a monodromy operator $N$,
following the work of Steenbrink \cite{st} and
Bloch-Gillet-Soul\'{e} \cite{bgs}. We need to use this complex to
define the Deligne cohomology in the case of strict semistable
reduction. To define it we need some preliminaries.

Let $Y=\coprod_{i=1}^{t} Y_i$ be the special fibre of dim $n$ with
$Y_i$ its irreducible components. For $I \subset \{1,2,...t\}$,
define
$$Y_I= \cap_{i \in I} Y_i$$
Let $r=|I|$ denote the cardinality of $I$. Define
$$Y^{(r)}:=\begin{cases} \XX & \text{ if } r=0 \\ \coprod_{|I|=r} Y_{I}& \text{ if } 1 \leq r \leq n \\
\emptyset & \text { if } r>n \end{cases}$$
For $u$ and $t$ with $1 \leq u \leq t < r$ define the map
$$\delta(u):Y^{(t+1)} \rightarrow Y^{(t)}$$
as follows. Let $I=(i_1,i_2,...i_t)$ with $i_1<i_2<...<i_t$. Let
$J=I-\{i_u\}$. There is an embedding
$$Y_I \rightarrow Y_J$$
and this induce maps on the cohomology and homology of these
varieties. Let $\delta(u)^*$ and $\delta(u)_*$ denote the
corresponding maps. They further  induce  the Gysin and
restriction maps on the cohomology and homology as follows.

Define
$$\gamma:=\sum_{u=1}^{r} (-1)^{u-1} \delta(u)_*$$
and
$$ \rho:=\sum_{u=1}^{r} (-1)^{u-1} \delta(u)^*$$
These maps have the properties that
\begin{itemize}
\item $\gamma^2=0$ \item $\rho^2=0$ \item $\gamma \cdot \rho +
\rho \cdot \gamma =0$
\end{itemize}

Let $i,j,k \in \ZZ$. Define, following \cite{co}(3.1)
$$K^{i,j,k}:=\begin{cases}
CH^{\frac{i+j-2k+n}{2}}(Y^{(2k-i+1)})\otimes \Q & \text { if } k
\geq
max(0,i) \\
0 & \text{ otherwise}
\end{cases}
$$
and let
$$K^{i,j}= \oplus_{k} K^{i,j,k} \text { and } K^n=\oplus_{i+j=n}
K^{i,j}$$
The maps $\rho$ and $\gamma$ induce differentials
$$\partial':K^{i,j,k} \rightarrow K^{i+1,j+1,k+1} \hspace{1in}
\partial'(a)=\rho(a)$$
$$\partial'':K^{i,j,k} \rightarrow K^{i+1,j+1,k} \hspace{1in}
\partial''(a)= - \gamma(a)$$
Further define
$$N:K^{i,j,k} \rightarrow K^{i+2.j,k+1}(-1) \hspace{1in} N(a)=a$$
Let $\partial=\partial'+\partial''$ on $K^{i,j}$. From the
definition we have $[\partial,N]=0$ and $\partial^2=0$.

Let $Cone(N):K^* \rightarrow K^*$ be the complex $K^*\oplus
K^*[-1]$ with differential
$$D(a,b)=(\partial(a),N(a)-\partial(b))$$

Consani \cite{co}[Prop 3.4] shows that this cone complex is
quasi-isomorphic to a complex of Chow groups of the fibre:

\begin{prop}[Consani]
Let $*$ be a fixed integer. The complex, for $q \in \ZZ$,
$$Cone(N:K^{q-2*,q-n} \rightarrow K^{q-2*+2,q-n})$$
is quasi-isomorphic to the following complex
$$\CC^{q}(*)=\begin{cases} CH^{q-*}(Y^{(2*-q)}) & \text {
if $q \leq *-1$}\\
CH^*(Y^{(q-2*)},\QL) & \text{ if $q \geq *$}
\end{cases}
$$
The differential $d_{\CC}$ is given by
$$d_{\CC}(a)=\begin{cases} d''(a) & \text{ if } q<*-1\\
-i^{*}i_{*}(a)& \text { if } q=*-1\\
d'(a) & \text{ if } q \geq *
\end{cases}
$$
\end{prop}

\section{ The `Deligne cohomology' at a finite prime}
Assume now that the residue field $k(v)$ is finite \footnote{This
is necessary to define the $L$-series and is not strictly
necessary at this point.}. We define the {\bf $v$-adic Deligne
Cohomology group} to be
$$H_{\D}^{q}(X_{/v},\Q(q-a)):=\begin{cases} CH^{q-a-1}(Y,q-2a-1)\otimes
\Q & \text { if $q-2a>1$} \\
\frac{Ker(i^*i_*:CH_{n-a}(Y^{(1)}) \rightarrow
CH_{n-a}(Y^{(1)}))}{Im(\gamma:CH_{n-a}(Y^{(2)}) \rightarrow
CH_{n-a}(Y^{(1)}))}\otimes \Q & \text{ if $q-2a=1$} \end{cases}
$$
Here $n$ is the dimension of $Y$. This is a $\Q$ vector space
which we will show has the expected properties assuming certain
conjectures. Note that if  $Y$ is non-singular and $q-2a>1$ the
Parshin-Soul\'{e} conjecture asserts that the higher Chow group is
finite, hence this space is $0$. When $q-2a=1$  the group is
$CH^a(Y) \otimes \Q$.

\begin{rem}
Consani shows that if one uses a certain complex of differential
forms to define $K^{i,j,k}$ and performs the same calculations,
one ends up with the Real Deligne cohomology as the graded pieces.
\end{rem}

\section{Properties of the Deligne Cohomology}

\subsection{Dimension:} The usual Real Deligne cohomology has
the property that its dimension is the order of the pole  of the
Archimedean    part of the $L$-function at a certain point on the
left of the critical point. Here we have a similar property. Let
$F^*$ be the geometric Frobenius and $N(v)$ the number of elements
of $k(v)$. The local $L$-factor of the $(q-1)^{st}$ cohomology
group is then
$$L_v(X,s)=(det(\textsl{I}-F^*N(v)^{-s}|H^{q-1}(\bar{X},\QL)^{I}))^{-1}$$
\begin{thm}[Consani] Let $v$ be a place of semistable reduction.
Assuming the weight-monodromy conjecture, the Tate conjecture for
the components  and the injectivity of the cycle class map on the
components $Y_I$, Parshin-Soul\'{e} conjecture and that $F^*$ acts
semisimply on $H^*(\bar{X},\QL)^{I}$. we have
$$dim_{\Q}(H_{\D}^{q}(X_{/v},\Q(q-2a)))=-ord_{s=a} L_{v}(X,s):=d_v$$
\end{thm}

\begin{proof} \cite{co}, Thm 3.5
\end{proof}

\begin{rem} Since the $L$-factor at a prime of good reduction does not have a pole at
$s=a$ when $q-2a>1$, the Parshin-Soul\'{e} conjecture can be
interpreted as the statement that the $v$-adic Deligne cohomology
has the correct dimension, namely $0$, even at a prime of good
reduction.
\end{rem}

\begin{rem} In the function fields setting and more recently, in
the setting of  p-adically uniformized varieties, the weight-monodromy
conjecture  is  a  theorem  \cite{de3}, \cite{it}.   Further,  in  the
p-adically uniformized case, the  variety is totally degenerate so the
Tate   conjecture  for   the  components   is  trivial,   so  assuming
semi-simplicity of the action of  the Frobenius and injectivity of the
cycle class  map on  the components, Consani's  theorem holds  in this
case.
\end{rem}

\subsection{Regulator maps}

To define the regulator map we use the localization sequence \cite{bl1}. It
is as follows. If $X, \XX$ and $Y$ are as before, and $q,a \in
\ZZ, q-2a>0 $ we have
$$\cdots \rightarrow CH^{q-a}(\XX,q-2a) \rightarrow CH^{q-a}(X,q-2a)
\stackrel{\partial}{\rightarrow} CH^{q-a-1}(Y,q-2a-1)$$
$$\rightarrow CH^{q-a}(\XX,q-2a-1) \rightarrow CH^{q-a}(X,q-2a-1)
\rightarrow \cdots$$
The usual regulator map should appear as the boundary map in the
`arithmetic' localization sequence.
$$ \cdots \rightarrow \widehat{CH}^{q-a}(\XX,q-2a) \rightarrow
CH^{q-a}(\XX,q-2a) \stackrel { r_{\D}}{\longrightarrow}
H^{2a+1}_{\D}(X_{/\R},\R(a+1))\rightarrow \cdots$$
As far as we are aware, higher arithmetic Chow groups have not
been defined in general, but this is known in the case when
$q-2a=1$.

In the finite prime case, there are two cases that have to be
considered.

\subsubsection{Case 1: $q-2a>1$}

We define the {\bf $v$-adic regulator map} to be the map
$\partial$.
$$r_{\D,v}:CH^{q-a}(X,q-2a) \stackrel{\partial}{\longrightarrow} H^{q}_{\D}(X_{/v},\Q(q-a))$$
In analogy with the Beilinson conjectures, we have the following
conjecture \vspace{.25in}

\noindent {\bf CONJECTURE A1:} The image
$$Im(r_{\D,v}(CH^{q-a}(X,q-2a))) \subset  H^{q}_{\D}(X_{/v},\Q(q-a))$$
is a full sub-lattice.

\subsubsection {Case 2: $q-2a=1$} In this case the Chow groups of the
components of the special fibre are not torsion, so the conjecture
has to be slightly modified. However, we get a conjecture which is
non trivial even in the case of good reduction.

One has a presentation of the Chow group
$$CH^*(Y)=Coker(\gamma:CH_{n-*}(Y^{(2)}) \rightarrow
CH_{n-*}(Y^{(1)}))$$
Recall that $i:Y \hookrightarrow \XX$ is the inclusion map of the
special fibre into the model. From the exactness of the
localization sequence, the image of $\partial$ lies in the kernel
of
$$i_*:CH_{n-*}(Y) \rightarrow CH_{n-*}(\XX)=CH^{*-1}(\XX)$$
In particular, it lies in the kernel of $i^*i_*:CH_{n-*}(Y^{(1)})
\rightarrow CH_{n-*}(Y^{(1)})$.

We define the {\bf $v$-adic regulator map} as before as  map
$\partial$
$$r_{\D,v}:CH^{a+1}(X,1) \stackrel{ \partial}{\longrightarrow} H^{2a+1}_{\D}(X_{/v},\Q(a+1))$$
and we have the following conjecture.

\vspace{.25in}

\noindent {\bf CONJECTURE A2:} The map
$$Im(r_{\D,v}(CH^{a+1}(X,1))) \subset  H^{2a+1}_{\D}(X_{/v},\Q(a+1))$$
is a full sublattice.

\subsubsection{The map $z^a_v$}

We have another map
$$z^a_v:CH^a(X) \longrightarrow H^{2a+1}_{\D}(X_{/v},\Q(a+1))$$
defined as follows. The map $j^*:CH^a(\XX) \rightarrow CH^a(X)$ is
surjective and we have the map $i^*:CH^a(\XX) \rightarrow
CH^a(Y)$. This induces a map $i^* \circ (j^*)^{-1}:CH^a(X)
\rightarrow CH^a(Y)$ which is well defined up to an element of the
image of $i^*i_*:CH^{a-1}(Y) \rightarrow CH^{a}(Y)$. So we have a
well defined map
$$\xi^a_v:CH^{a}(X) \rightarrow \frac{Ker(\rho:CH^a(Y^{(1)})
\rightarrow CH^{a}(Y^{(2)}))}{Im(i^*i_*)} \otimes \Q$$
There is always a morphism $\tau$ from the group on the right to
the Deligne cohomology induced by the cap product $\cap[Y]$
\cite{bgs}, pg. 454., and we define
$$z^a_v:=\tau \circ \xi^a_v$$
In some  instances, for example, if $Y$ is smooth \cite{bgs}[Prop
2.3.3], the map $\tau$ is known to be an isomorphism. In several
other instances \cite{bgs}[Section 6], the map is known to be an
isomorphism after going to cohomology.

\subsubsection{ The Archimedean case}

The Archimedean case of {\bf Conjecture A} is Beilinson's
Hodge-$\D$-conjecture. However, this is known to be false in
general \cite{sms}. It was suggested by N.Fakhruddin \cite{naf}
that perhaps similar methods can be used to show that {\bf
Conjecture A} is false over p-adic fields, but may still  hold
over dvr's whose quotient field is contained in the algebraic
closure of a number field. As in all the examples we have the
conjecture is true, we have left it as it is. However, in all our
examples the Hodge-$\D$-conjecture is also known.

The Archimedean version of the map $z^a$ comes from the usual
cycle class map to Betti cohomology  and the long exact sequence
relating the Betti, de Rham and Deligne cohomologies.

\section{The $S$-integral Beilinson conjecture}
The conjectures above can be combined with the Hasse-Weil
conjecture to formulate an $S$-integral version of the Beilinson
conjectures. Let $X$ be a variety defined over a number field $K$
which has at worst semi-stable reduction at all places ( this
restriction may not be so important but at the moment it is
necessary ). Let $\OO_K$ be the ring of integers of $K$ and let
$\bfinfty$ denote the set of Archimedean    places. Let
$L_{v}(X,s)$ be a cohomological $L$-function at the place $v$ for
the $(q-1)^{st}$ cohomology group as defined  above, where for $v$
Archimedean    it is defined as in \cite{rss} pg 4. Define the
completed $L$-function $\Lambda(X,s)$ as follows:
$$\Lambda(X,s)=A^{s/2}\prod_{v} L_v(X,s)$$
where $A$ is a generalized conductor.

The {\bf Standard Conjectures} of Grothendieck, generalizing those
of Hasse-Weil and Serre,  asserts that this function has a
meromorphic continuation to the entire complex plane
$$\Lambda(X,s)=\pm \Lambda(X,q-s)$$
For a set of places $S$ containing the Archimedean    places we
define the $S$-integral $L$-function as follows:
$$L_S(X,S)=\prod_{s \notin S} L(X,s)$$
so
$$\Lambda(X,s)=L_S(X,s) \prod_{v \in S} L_v(X,s)$$
The the above conjectures can be combined with the usual Beilinson
conjectures to give the following $S$-integral Beilinson
conjectures: Let $c_{\bfinfty}(X,q,a)$ be the number appearing the
usual Beilinson conjectures, coming from isomorphism of the two
$\Q$-structures and for a finite set of places $S$ containing the
Archimedean places, let $\XX_S$ be a model over $\OO_K[S]$ with at
worst semistable reduction.

\vspace{.25in}

\noindent{\bf CONJECTURE B1:} Let $S$ be a finite set of places
containing all the Archimedean    places and $q-2a>1$.Then

\begin{itemize} \item $ord_{s=a} L_S(X,s)=dim_{\Q}
CH^{q-a}(\XX_{S},q-2a) \otimes \Q$

\item The map
$$\oplus_{v \in S \backslash \bfinfty}
r_{\D,v} \otimes \Q  \bigoplus \oplus_{v|\bfinfty}
r_{\D,v} \otimes \R:CH^{q-a}(\XX_S,q-2a)
\longrightarrow \bigoplus_{v
\in S} H^{q}_{\D}(X_{/v},\Q(q-a))$$
is an isomorphism.

\item $L_S^*(X,a)\sim_{\Q^*} c_{\bfinfty}(X,q,a)\prod_{v \in S
\backslash \bfinfty} (\log(N(v)))^{d_v} $

\end{itemize}

\vspace{.25in}

\noindent {\bf CONJECTURE B2:} Let $S$ be a finite set of places
containing the Archimedean    places and $q-2a=1$. Let $B^a(X)$
denote the group $(CH^a(X)/CH^a_{hom}(X))$ and $z^a$ the cycle
class map to Deligne cohomology induced by the cycle class map to
$H^{2a}_{B}(X,\R(a))^{(-1)^a}$

Then
\begin{itemize}

\item $ord_{s=a} L_S(X,s)=dim_{\Q} CH^{q-a}(\XX_{S},1)\otimes \Q$

\item $ord_{s=a+1} L_S(X,s)=-dim_{\Q} B^a(X)$ ( Tate's Conjecture)

\item The map
$$\oplus_{v \in S \backslash \bfinfty } r_{\D,v}\otimes \Q \oplus_{v| \infty} r_{\D,v}\otimes \R \oplus z^a
:CH^{q-a}(\XX_{S},1) \oplus B^a(X)\otimes \Q \rightarrow \oplus_{v
\in S} H^q_{\D}(X/v,\Q(q-a))$$
is an isomorphism.

\item $L_S^*(X,a)\sim_{\Q^*} c_{\bfinfty}(X,q,a) \prod_{v \in
S\backslash \bfinfty} (\log(N(v)))^{d_v}$
\end{itemize}
\section{Function Fields}

In the case of function fields, we can formulate a conjecture
along the lines of  {\bf Conjecture B} using the fact that the
primes at ${\bfinfty}$ are just finite primes.

Let $\F_q$ be the finite field of order $q$ and $K$ a function
field over $\F_q$ of transcendence degree $1$. Let $S$ be a non
empty finite set of primes of $K$. Then there exists an element
$x$ in $K$ whose poles are precisely the elements of $S$. Let
$\OO_S$ denote the integral closure of $\F_q[x]$ in $K$. The prime
$\infty$ is the prime $1/x$ of $\F_q[x]$ and the set $S$ is the
set of primes lying over $\infty$. Let $\XX_S$ be a model of $X$
over $\OO_S$ with at worst semi-stable reduction.

\vspace{.25in}

\noindent{\bf CONJECTURE B1FF:} Let $S$ be a finite set of places
and $\XX_S$ as above and $q-2a>1$ Then
\begin{itemize}
\item $ord_{s=a} L_S(X,s)=dim_{\Q} CH^{q-a}(\XX_{S},q-2a) \otimes
\Q$

\item The map
$$\bigoplus_{v \in S} r_{\D,v}\otimes \Q:CH^{q-a}(\XX_S,q-2a) \rightarrow \bigoplus_{v \in S} H^{q}_{\D}(X_{/v},\Q(q-a))$$
is an isomorphism.

\item $L_S^*(X,a)\sim_{\Q^*} \prod_{v \in S} (\log(N(v)))^{d_v} $

\end{itemize}
\vspace{.25in}

In the case $q-2a=1$, let $B^a(X)$ denote the group
$(CH^a(X)/CH^a_{hom}(X))$ and let $z^a_S$ denote the map
$$z^a_S:=\oplus_{s\in S} z^a_v:B^a(X) \rightarrow \oplus_{s\in S}
H^{2a+1}_{\D}(X/v,\Q(a+1))$$
induced by the restriction maps $CH^a(X) \rightarrow CH^a(Y_v)$.

\vspace{.25in}

\noindent {\bf CONJECTURE B2FF:} If $q-2a=1$, Let $S$ be a finite
set of places and $\XX_S$ as above. Then
\begin{itemize}

\item $ord_{s=a} L_S(X,s)=dim_{\Q} CH^{q-a}(\XX_{S},1)\otimes \Q$

\item $ord_{s=a+1} L_S(X,s)=-dim_{\Q} B^a(X)\otimes \Q$ ( Tate's
Conjecture )

\item The map
$$\bigoplus_{v \in S} r_{\D,v} \otimes \Q \oplus z^a_S\otimes \Q
:CH^{a+1}(\XX_{S},1) \oplus B^a(X) \rightarrow
\bigoplus_{v \in S} H^{2a+1}_{\D}(X/v,\Q(a+1))$$
is an isomorphism.

\item $L_S^*(X,a)\sim_{\Q^*} \log(q)^{dim(CH^{a+1}(\XX_S,1))}$
\end{itemize}
Notice that the only transcendental term is a power of $\log(q)$
which comes from the residue of the local $L$-function.

\section{A Special value conjecture}

The considerations of the previous sections allow us to formulate
a special value conjecture in the function field case. Perhaps
this is the same as the Bloch-Kato conjecture \cite{bk} but we
have not verified that. Formulating the conjecture here is a little easier as in the number field case, Beilinson's regulator map is only defined up to $\Q^*$.

\subsection{$\ZZ$-structures}

To conjecture an expression for the special value, we use the fact
that the Deligne cohomology at a finite prime has two $\ZZ$
structures, one natural and the other conjectural. The first $\ZZ$
structure comes from the integral Chow group
$CH^{q-a-1}(Y_v,q-2a-1)$. {\bf Conjecture A} asserts that the
image of the regulator map gives a second $\ZZ$ structure, which
is a subgroup of the integral Chow group.

From the earlier conjectures when $q-2a>1$, the map
$$CH^{q-a}(X,q-2a)\otimes \Q \stackrel{\oplus r_{\D,v}}{\longrightarrow} \bigoplus_{v}
CH^{q-a-1}(Y_{v},q-2a-1)\otimes \Q$$
is an isomorphism. Similarly, when $q-2a=1$, the map
$$CH^{a+1}(X,1)\otimes \Q \oplus B^a(X)\otimes \Q \stackrel{\oplus r_{\D,v} \oplus
z^a_v}{\longrightarrow} \bigoplus_{v} CH^{a}(Y_v) \otimes \Q$$
is an isomorphism. However, the maps are  defined integrally, so
the kernel and cokernel are torsion. We conjecture that they are
actually finite. Let $b(X,q,a)$ denote the order of the kernel and
$c(X,q,a)$ denote the order of the cokernel.

\subsection{ A special value conjecture}
We then have the following conjecture:
 \vspace{.25 in}

\noindent{\bf CONJECTURE CFF:}
\begin{itemize}
\item C1: If $q-2a>1$ then
$$ord_{s=a} \Lambda(X,s)=dim_{\Q} CH^{q-a}(X,q-2a)\otimes \Q$$

\item C2: If $q-2a=1$, the space $CH^{q-a}(X,1)\otimes \Q$ is {\bf
not} finite dimensional, but
$$ord_{s=a} \Lambda(X,s)=-dim_{\Q} B^a(X)\otimes \Q$$
\end{itemize}
Further,
$$\Lambda^*(X,a)=\pm \frac{c(X,q,a)}{b(X,q,a)}\log(q)^{ord_{s=a} \Lambda(X,s)}$$
\begin{rem}
One can also formulate an exact value conjecture in the
$S$-integral case and this can be viewed as the conjecture when
$S$ is the set of all primes.
\end{rem}

\section{ Examples}

\subsection {$K_1$ of fields}

This is the only case where all the conjectures are known. Here
$q=1$ and $a=0$. Let $K$ be a number field and $\Lambda$ the
completion of $\OO_K$ at a finite place $v$. Let $\eta$ denote the
generic point of $Spec(\Lambda)$.

\subsubsection{Local Case - Conjecture A:} As $q-2a=1$, we are in case $2$.
Here the $L$-function is
$$L_v(Spec(K),s):=L_{v}(H^0(Spec(\Lambda),\QL(0)),s)=\frac{1}{1-N(v)^{-s}}$$
which has a pole at $s=0$. The conjecture A2 then asserts that the
rank of the image of the regulator map is $1$.

The localization sequence gives
$$\cdots \rightarrow CH^1(Spec(\Lambda),1) \rightarrow CH^1(\eta,1)
\rightarrow CH^0(v) \rightarrow CH^1(Spec(\Lambda)) \rightarrow
CH^1(\eta)\rightarrow 0$$
Since $\Lambda$ is a principal ideal domain, $CH^1(\Lambda)=0$, so
the regulator map is surjective and the image has of the regulator
map is  rank $1$. So this is a case when the conjecture is true.
Notice this is a case in which the variety has good reduction.

The statement that the map
$$\R^* \stackrel{\log||}{\longrightarrow} \R$$
is surjective, and similarly for $\C^*$, can be viewed as an {\em
Archimedean    local case}.

\subsubsection{Global Case - Conjectures B and C :} Let
$$\zeta_{K,S}(s)=\prod_{v \notin S} L_v(Spec(K),s)$$

{\bf Conjecture B2} is the usual $S$-unit theorem of Dirichlet.
The fact that the regulator map is surjective follows from the
finiteness of the class number.

In this case one also has an expression for the special value
along the lines of {\bf Conjecture CFF} which follows from the
class number formula and the functional equation. Here
$$L_{\infty}(Spec(K),s)=2^{- r_2s}\pi^{-\frac{ns}{2}}\Gamma(s/2)^{r_1}\Gamma(s)^{r_2}$$
and
$$\Lambda_K(s)=|D(K)|^{s/2} L_{\infty}(Spec(K),s)\zeta_K^{\bfinfty}(s)$$
where $D(K)$ is the discriminant. $\Lambda_K(s)$ satisfies the
functional equation
$$\Lambda_K(s)=\Lambda_K(1-s)$$
and has a simple pole at $s=0$ with residue
$$\Lambda^*_K(0)=Res_{s=0}\Lambda_K(s)=-\frac{R_K h_K}{w_{K}}$$
Here $R_K$ is the classical regulator, $h_K$ is the class number
and  $w_K$ is the number of roots of unity. $h_K$ is the order of
the cokernel of the non-archimedean regulator map and $w_K$ is the
order of the kernel, so the formula is very similar to that
expected in {\bf Conjecture CFF}. In this case there is a way of
defining a second integral structure on the Deligne cohomology and
$R_K$ is the cokernel of that map.

The simple pole at $s=0$ shows that {\bf Conjecture C2} holds in
this case as well.

\subsubsection{Function field case}

In this case the precise analogue of all the conjectures is known
and can be found in the book of Rosen \cite{ro}. If $K=\F_q(T)$,
Rosen (pg 244) shows, for example, that the completed zeta
function $\Lambda_{K}(s)$ has a simple pole at $s=0$  and
$$\Lambda_{K}^*(0)=-\frac{h_K}{(q-1)\log(q)}$$
which is precisely what is expected by {\bf Conjecture C2} as
$h_K$ is the class number and is the order of the cokernel, while
the kernel is the number of elements of finite order in $K^*$
which is $q-1$, the cardinality of $\F_q^*$.

\subsection{ $K_2$ of curves }

Here we consider the case when $q=2$ and $a=0$  so $q-2a=2$.

\subsubsection{Local Case - Conjecture A} Here at a prime of good
reduction the Deligne cohomology group is trivial, so that case is
true for trivial reasons. At a prime of semistable reduction
Ramakrishnan  \cite{ra} following Bloch and Grayson \cite{bg} for
elliptic curves  defined a regulator map which coincides with
ours. Bloch and Grayson showed that Conjecture A is true for a
prime of split semistable reduction of an elliptic curve over
$\Q$. A more precise computation  of the $v$-adic regulator at
such a prime can be found in the work of Schappacher and Scholl
\cite{ss}.

\subsubsection {Global Case: Conjecture B} In general not much is know for
curves. The usual Beilinson conjecture is known for some types of
elliptic curves over $\Q$ and for modular curves. The $S$-integral
version was originally formulated by Ramakrishnan, but still
remains to be proved.

\subsubsection{ Function field case }

The function field we consider is $K=\F_q(T)$. This case when
$S=\{1/T\}$ is done in the works of Kondo \cite{ko} and P\`{a}l
\cite{pal}. They define and compute the regulator map for $K_2$ of
Drinfel'd modular curves and, via the modular parametrization for
the product of two elliptic curves.

\subsection{ $K_1$ of products of  modular curves and elliptic curves }

This is the case when $q=3$ and $a=1$. In this case one expects a
contribution from the primes of good reduction as well. From the
modularity theorem, what is proved here can be used to prove a
statement about products of elliptic curves over $\Q$.

\subsubsection{ Local Case - Conjecture A}

Let $E$ and $E'$ be two elliptic curves over $\Q$ and $v$ a finite
place. We are interested $CH^2(X,1)$ where $X=E \times E'$. Since
they are modular, we can compute the $L$-function of $H^2(X)$.
There are several cases to be considered.

First, assume $E$ and $E'$ are not isogenous. Suppose both have
good reduction at $v$ and $E_v$ and $E'_v$ are not isogenous. Then
$H^3_{D}(X/v,\Q(2))$ is two dimensional. The map $r_{\D,v}$ is
surjective as we have a map $CH^1(X,1) \otimes CH^1(X,0)
\rightarrow CH^2(X,1)$. The elements in the image are called {\em
decomposable}. Consider the cycles $(E \times pt)\otimes 1)$ and
$(pt \times E')\otimes 1)$. These map surjectively onto the
Deligne cohomology so that case is taken care off. Suppose $E_v$
and $E'_v$ are isogenous. Then $H^3_{D}(X/v,\Q(2))$ is $4$
dimensional - the new cycles coming from the graph of the isogeny
and the graph of the isogeny composed with the Frobenious. Spie\ss
\cite{sp} shows how to find elements in $CH^2(X,1)$ whose
regulator is these extra cycles.

If $E$ and $E'$ are isogenous, then $E_v$ and $E'_v$ are
necessarily isogenous and so we are in case $2$ above. It should
be remarked that over larger residue fields, the rank of the
Deligne cohomology could be as large as $6$ owing in the case of
supersingular reduction. However, if $E$ is an elliptic curve over
$\Q$ and $E_p$ is supersingular, the extra endomorphisms are not
defined over $\F_p$.

If $E_v$ has good reduction and $E'v$ does not, the $L$-function
computation suggests that the rank of the Deligne cohomology is
$2$ and hence we can use the decomposable cycles to prove
surjectivity. If both $E_v$ and $E'_v$ have semi-stable reduction,
then $L$-function shows that the rank of the Deligne cohomology is
$3$. If $E$ and $E'$ are not isogenous, two of the elements can be
obtained from the decomposable elements. The third is constructed
using the modular parametrization, a construction similar to that
of Mildenhall \cite{mi}. Finally, if $E$ and $E'$ are isogenous,
and $v$ is a prime of semistable reduction, the Deligne cohomology
is once again 3 dimensional and the surjectivity of the regulator
map comes from the 3 decomposable elements.

The achimedean case, namely the Hodge-$\D$-conjecture is due to
Lewis and Chen \cite{cl}.

\subsubsection{ Global Case - Conjecture B:}

When $S=\bfinfty$ and $X$ is the product of two modular curves,
Beilinson \cite{be} proved that the regulator map is surjective
and showed that the value of the $L$-function is equal to
$c_{\bfinfty}(X,3,1)$ up to a rational number.

From the work of Mildenhall, one can deduce the weak $S$-integral
conjecture for any finite set of primes $S$ of good or semi-stable
reduction.

The exact value conjecture is still unresolved, in the case of $E$
and $E'$ isogenous, the special value is related to the symmetric
square of $L$-function of a modular form and there is a lot of
work on this due to Flach \cite{fl} and others. In \cite{ba-sr},
we obtain an exact formula in the non-isogenous case  which may
have some bearing on the conjecture.

\subsubsection{Function Field case}

The case of conjecture $B$ with $S=\{1/T\}$ for $K_1$ of the
product of Drinfeld modular curves was carried out by us
\cite{cs}. We show, for example, that the special value of the
$L$-function is $\log(q)$ up to an element of $\Q^*$, as predicted
by the conjecture, though we actually have a more precise
expression for the special value.

\vspace{.2in}

\vskip .1in

\noindent {\bf Ramesh Sreekantan}, Department of Mathematics,
University of Toronto, Canada.

\noindent email: ramesh\@@math.toronto.edu

\end{document}